\RequirePackage{amsmath}
\documentclass{svjour3}                     
\smartqed  
\usepackage{graphicx}
\usepackage{latexsym}

\usepackage{amsthm}
\usepackage{amsmath}
\usepackage{amssymb}
\usepackage{breqn}
\usepackage{appendix}
\usepackage{enumerate}
\newtheorem{thm}{Theorem}[section]
\newtheorem*{thm*}{Theorem}

\newtheorem{lem}[thm]{Lemma}
\newtheorem*{lem*}{Lemma}
\newtheorem{cor}[thm]{Corollary}

\newtheorem*{prop*}{Proposition}

\newcommand{\R}{\mathbb{R}}

\newcommand{\N}{\mathbb{N}}

\newcommand{\T}{\mathcal{T}}
\newcommand{\Leray}{\mathcal{P}}

\begin{document}

\title{Global regularity of the 2D HVBK equations}

\author{Pranava Chaitanya Jayanti \and Konstantina Trivisa}

\institute{P. C. Jayanti \at
              Department of Physics / Centre for Scientific Computation And Mathematical Modeling \\
              University of Maryland, College Park, USA\\
              \email{jayantip@cscamm.umd.edu}\\
              ORCID iD: https://orcid.org/0000-0001-8151-406X \\
           \and
           K. Trivisa \at
              Department of Mathematics / Institute for Physical Science and Technology / Centre for Scientific Computation And Mathematical Modeling \\
              University of Maryland, College Park, USA }

\date{Received: date / Accepted: date}

\maketitle

\begin{abstract}
The Hall-Vinen-Bekharevich-Khalatnikov (HVBK) equations are a macroscopic model of superfluidity at non-zero temperatures. For smooth, compactly supported data, we prove the global well-posedness of strong solutions to these equations in $\R^2$, in the incompressible and isothermal case. The proof utilises a contraction mapping argument to establish local well-posedness for high-regularity data, following which we demonstrate global regularity using an analogue of the Beale-Kato-Majda criterion in this context. In the appendix, we address the sufficient conditions on a 2D vorticity field, in order to have a finite kinetic energy.

\keywords{HVBK equations \and Superfluids \and Global well-posedness \and Navier-Stokes equations}
\end{abstract}

\pagebreak

\section{Introduction}  \label{intro}

Most substances upon isobaric cooling transition from a gas to liquid, before eventually turning into a solid phase (or in some cases, a variety of them). Helium is an exception - at pressures below 25 bars, liquid Helium-4 transforms into a \textit{superfluid} phase when cooled across the \textit{lambda line} (approximately 2.17 Kelvin)\footnote{Helium-3 also displays a superfluid phase, albeit at a significantly lower temperature ($\sim$2 milliKelvin) due to its fermionic nature. Most experimental research has focused on Helium-4.}. This superfluid phase was experimentally discovered \cite{Kapitza1938Viscosity-Point,Allen1938FlowII} over 80 years ago and has since been an important subject of interest to, and investigation by, the physics community. As far as a theoretical explanation of superfluidity goes, there are several pieces to the puzzle that work in limited ranges of validity; however, a single universal theory that explains everything to reasonable satisfaction continues to elude us. The most well-known models currently in use may be classified in many ways \cite{Barenghi2014IntroductionTurbulence}, the one on the basis of length scales being of interest here.

One of the first theories was by Fr\"ohlich (in 1937), proposing a model of order-disorder transition where a fraction of the helium atoms are trapped in a lattice that denotes the ground state, while the remaining atoms are in an excited state. In 1938, London \cite{London1938TheDegeneracy} modified this theory by suggesting that the ground state could actually be a degenerate Bose gas. Shortly thereafter, Tisza \cite{Tisza1938TransportII} worked on the details of London's model, suggesting modifications to match with experimental data. In 1941, Landau \cite{Landau1941TheoryII} presented his famous two-fluid model, which would later earn him the Nobel Prize. This was a semi-classical model in which he considered the superfluid to be described by a classical velocity field, and the normal fluid to be the excitations (phonons and rotons) of this underlying superfluid field. The phonons were the excitations corresponding to linear modes (quantized sound waves), while the rotons were those of rotational modes that were identified with a local minimum in the energy-momentum diagram. He emphasised that the two fluids cannot really be compared to a (classical) multiphase flow where each point in spacetime can be \textit{uniquely} identified with a given phase. Landau's model made experimentally testable and qualitatively accurate predictions, particularly with regards to the critical velocity $-$ when the superfluid is moving below the critical velocity, it is dissipation-free; if it moves faster, then it has enough energy to excite some phonons and rotons, which can dissipate energy via interactions with the walls of the container.

These works, along with the identification of the superfluid phase with a BEC, led to the increasing use of the Gross-Pitaevskii equation (GPE), a nonlinear dispersive PDE. The wavefunction that solves this equation is an order parameter which vanishes at superfluid vortices\footnote{The existence of these vortices was predicted by Onsager \cite{Onsager1953IntroductoryTalk} and Feynman \cite{Feynman1955ApplicationHelium} in the mid-50's, and were directly observed in 2006 \cite{Bewley2006VisualisationVortices} by pinning tracer particles to the vortices to visualise them.} (also known as ``topological defects") \cite{Paoletti2011QuantumTurbulence}. The main limitations \cite{Roberts2001TheSuperfluidity} of the GPE are that it is valid only for low-energy scattering of the condensate particles, and also only at absolute zero. The low-energy limit leads to a local potential (Dirac delta), and must therefore be replaced with more general, non-local potentials to describe the bosonic interaction dynamics at higher energies. At $T=0$K, the helium is completely condensed into the ground state and is composed only of the inviscid superfluid phase. As previously mentioned, this phase has no viscosity, and at sufficiently low velocities, no excitations either. The nonlinear Schr\"odinger equation (NLS), as the GPE is known in mathematics literature, has been studied extensively by mathematicians for well-posedness \cite{CollianderWell-posednessEquations} and scattering (see \cite{Tao2006NonlinearAnalysis,Dodson2016Global2} and references therein) results over the last few decades. However, as the system is heated to non-zero temperatures, the presence of (and interactions between) the two fluids must be accounted for. This has led to the emergence of various models \cite{Berloff2014ModelingTemperatures}, each of which work well at certain characteristic length scales (micro-, meso- and macro-scales). The basic idea underlying all these models is that an interaction between the two fluids would cause momentum and energy to flow between them, and the normal fluid would continuously dissipate energy through viscosity.

\paragraph{Micro-scale models:} Micro-scale models describe the system at a quantum level (length scales $\ll$ inter-vortex distance), and are essentially modifications of the NLS. One way is to simply ignore the normal fluid altogether and represent its dissipative effect on the superfluid by imposing a high-frequency cutoff \cite{Kobayashi2005KolmogorovDissipation}, or adding nonlinear damping terms \cite{Antonelli2010GlobalDamping}, or fractional diffusion operators \cite{Darwich2016SomeDissipation,Darwich2018GlobalDissipation}. Another approach is to assume a given velocity fluid for the normal fluid and include its effect on the superfluid in the form of a material transport term (material derivative) of the density of the superfluid \cite{Carlson1996AVortices}. Finally, in this hierarchy of models, the most sophisticated may be one by Pitaevskii \cite{Pitaevskii1959PhenomenologicalPoint}, in which a set of equations coupling (bi-directionally) the superfluid and the normal fluid were derived from first principles.

\paragraph{Meso-scale models:} At non-zero temperatures, the superfluid contains some excitations (vortices) that can interact with the normal fluid. To describe the dynamics at length scales comparable to the inter-vortex spacing, Schwarz \cite{Schwarz1978TurbulenceCounterflow,Schwarz1985Three-dimensionalInteractions,Schwarz1988Three-dimensionalTurbulence} utilised a vortex model where he treated the superfluid to be incompressible and irrotational, except at the vortices. Each vortex is subject to a drag force from the superfluid velocity field (of the remaining vortices) and the normal fluid velocity field. Apart from this, numerical simulations of the system are designed to include \textit{vortex reconnections}, predicted to occur in superfluids by Feynman \cite{Feynman1955ApplicationHelium}. Vortex reconnection also plays an important role in quantum turbulence \cite{Vinen2006AnTurbulence,Paoletti2011QuantumTurbulence,Barenghi2014IntroductionTurbulence}, since reconnection leads to the formation of Kelvin waves along the vortices, which interact nonlinearly with each other resulting in dissipation through phonon emission.

\paragraph{Macro-scale models:} When the length scales of interest are much larger than the inter-vortex spacing, the superfluid can be modelled like a classical inviscid fluid. This prompted an effort to express the manifestation of quantum effects on large scales, discarding the discreteness of the vortex filaments for a continuum viewpoint, thus giving way to the Hall-Vinen-Bekharevich-Khalatnikov (HVBK) model \cite{Barenghi2001QuantizedTurbulence}. In this model, the viscous normal fluid is described using the Navier-Stokes equations, while the inviscid superfluid is represented by the Euler equations. This model is known to work well in situations when the vortex filaments are more or less aligned with each other, and is thus used to study superfluid Couette flow in rotationally symmetric domains \cite{Barenghi1988TheII,Henderson2004SuperfluidAnnulus,Peralta2009SuperfluidFlow}. In this paper, we prove the global regularity of strong solutions to the HVBK equations in two dimensions. \\

Although the original HVBK model treated the superfluid as inviscid, this leads to significant difficulties during numerical simulations since there would not be small-scale dissipation, forcing one to resolve extremely high frequencies. There have been multiple ways reported in the literature to get around this: suppressing high-frequencies in the superfluid velocity \cite{Salort2011MesoscaleTurbulence,Barenghi2014ExperimentalFluid}, using subgrid models \cite{Tchoufag2010EddyII} similar to the ones from Large Eddy Simulation (LES), and adding an artificial superfluid viscosity \cite{Roche2009QuantumCascade}. In this work, we follow the last method, i.e., by treating the superfluid as viscous. We will now state and describe the model in consideration. \\

\subsection{Notation}
We use the subscript $x$ on a Banach space to denote the Banach space is defined over $\R^2$. For instance, $L^p_x$ stands for Lebesgue space $L^p(\R^2)$, and similarly for the Sobolev spaces $H^m_x$. In the same way, spaces over time are denoted with a subscript $t$, like in $L^p_t$. Here, the time interval is understood from context, and will be specified otherwise. We also use subscripts on integrals to denote the variable being integrated; for example, $\int_x$ would stand for $\int_{\R^2}$, $\int_t$ for $\int_0^T$ and so on. \\

Since we are in 2D, the curl operator (acting on a 2D field) is often replaced by the \textit{perpendicular gradient} operator, denoted by: 
$$\nabla^{\perp} := \begin{bmatrix} 0 & 1\\ -1 & 0 \end{bmatrix}\cdot\nabla$$

We also use the notation $X\lesssim Y$ to imply that there exists a positive constant $C$ such that $X\le CY$. The constant may depend on various parameters ($k_1,k_2$, say), in which case it will be denoted using a subscript as $X\lesssim_{k_1,k_2} Y$ or $X\le C_{k_1,k_2}Y$.

\subsection{Organisation of the paper}

In Section \ref{mathematical model}, we present and discuss the mathematical model, along with statements of the main results. Section \ref{apriori estimates} contains some apriori estimates, followed by the proof of local well-posedness and pressure regularity in Section \ref{local well-posedness}. The global well-posedness for high-regularity data, including the analogue of the Beale-Kato-Majda condition, constitute Section \ref{global well-posedness}. Finally, we consider global regularity for smooth initial data in Section \ref{GWP for smooth data section}. The appendix contains a simple lemma, showing weaker and sufficient conditions on the vorticity in order to have a finite kinetic energy in 2D. \\

\begin{acknowledgements}
P.C.J. was supported by the Kulkarni Summer Research Fellowship at the University of Maryland during the summer of 2019, when parts of this work were completed. P.C.J. is also thankful to Jacob Bedrossian for bringing to the authors' attention the question addressed in the appendix. 
K.T. gratefully acknowledges the support by the National Science Foundation under the award DMS-1614964.
\end{acknowledgements}

\section{Mathematical model and main results} \label{mathematical model}
We first consider the general form of the incompressible HVBK equations as given in section 3.3 of \cite{Barenghi2001IntroductionTurbulence}. The equations \eqref{HVBK generalized} are repeated below for convenience (with slightly altered notation).

\begin{equation} \tag{HVBK} \label{HVBK generalized}
    \begin{gathered}
        \nabla\cdot u_n = 0 \quad , \quad \nabla\cdot u_s = 0 \\
        \partial_t u_n + (u_n\cdot \nabla)u_n = -\frac{1}{\rho} \nabla p_n - \frac{\rho_s}{\rho_n}S\nabla\Theta + \nu_n \Delta u_n + \frac{\rho_s}{\rho}F \\
        \partial_t u_s + (u_s\cdot \nabla)u_s = -\frac{1}{\rho} \nabla p_s + S\nabla\Theta + T - \frac{\rho_n}{\rho} F \\
        \omega_s = \nabla \times u_s \quad , \quad T = -\lambda_s \omega_s \times (\nabla \times \Hat{\omega}_s) \\
        \lambda_s = \frac{\Gamma}{4\pi} \ln{\left( \frac{b_0}{a_0} \right)} \quad , \quad b_0 = \left( \frac{\Gamma}{2\lvert\omega_s\rvert} \right)^{\frac{1}{2}} \\
        F = \frac{B}{2}\Hat{\omega}_s\times[\omega_s \times (u_n - u_s) + T] + \frac{B'}{2}[\omega_s \times (u_n - u_s) +T]
    \end{gathered}
\end{equation}

\medskip

In the above equations, subscripts $n$ and $s$ are used to denote the normal and superfluid components. $u_n$ and $u_s$ are the fluid velocities, $\omega_s$ is the superfluid vorticity,  $\Hat{\omega}_s$ is the unit vector along the vortex filament, $\rho_n$ and $\rho_s$ their (constant) densities, $\rho = \rho_n + \rho_s$ is the total density, $p$ is the pressure, $S$ is the entropy, $\Theta$ is the temperature, $\nu_n$ is the normal fluid viscosity, $T$ is the vortex tension term (force associated with bent vortex filaments), $F$ is the ``mutual friction" (the coupling between the two fluids), $\Gamma$ is the quantum of circulation in a superfluid vortex filament, $B$ and $B'$ are temperature-dependent positive friction coefficients, and $a_0$ is a positive constant that denotes the vortex core radius. \\

Now, we consider the isothermal case in 2D: the temperature gradient terms drop out from both equations; and the direction of the vorticity is fixed, so vortex lines do not bend, i.e., the vortex tension term is not relevant. As mentioned before, we will also be including a superfluid viscosity. Finally, we absorb $\rho$ in the definition of pressure and arrive at the following form of the equations.

\begin{gather}
    \nabla\cdot u_n = 0 \quad , \quad \nabla\cdot u_s = 0 \label{incompressibility} \\
    \partial_t u_n + (u_n\cdot \nabla)u_n = -\nabla p_n + \nu_n \Delta u_n + \frac{\rho_s}{\rho}F \label{normal fluid} \\
    \partial_t u_s + (u_s\cdot \nabla)u_s = -\nabla p_s + \nu_s \Delta u_s - \frac{\rho_n}{\rho} F \label{superfluid} \\
    F = -\frac{B}{2}\lvert\omega_s\rvert(u_n - u_s) + \frac{B'}{2}\omega_s\times (u_n - u_s) \label{mutual friction}
\end{gather}

\medskip

This is the system that will be dealt with in this paper, and is very similar to those considered in \cite{Roche2009QuantumCascade,Tchoufag2010EddyII,Verma2019TheModel}. Note that the first term in the mutual friction indeed acts like a (nonlinear) drag, and works to equalise the two velocity fields. The second term is a transport in the direction perpendicular to the relative velocity of the two fluids, so it is not a retarding force in itself. However, that is not concerning since it will be shown that this second term vanishes in the energy estimates, and behaves the same way as the first term as far as the contraction mapping argument goes. Furthermore, it is also evident that multiplying the normal fluid equation by $\rho_n$ and the superfluid equation by $\rho_s$, and adding the resulting equations leads to a cancellation of the mutual friction term. This is not surprising, as the friction force is internal to the entire system (of normal fluid and superfluid) and cannot alter the total momentum. \\

We seek mild solutions to \eqref{incompressibility}-\eqref{mutual friction}, expressed formally as:

\begin{gather}
    u_n(t) = e^{\nu_n t\Delta}u_n^0 - \int_0^t e^{\nu_n (t-\tau)\Delta}\Leray\left[u_n\cdot\nabla u_n - \frac{\rho_s}{\rho}F \right] d\tau \label{mild solution normal fluid} \\
    u_s(t) = e^{\nu_s t\Delta}u_s^0 - \int_0^t e^{\nu_s (t-\tau)\Delta}\Leray\left[u_s\cdot\nabla u_s + \frac{\rho_n}{\rho}F \right] d\tau \label{mild solution superfluid}
\end{gather}

Here, $\Leray$ is the Leray projector that gives as output the divergence-free part of a vector field; it is used to eliminate the pressure term. \\

It is a classical result that the incompressible Navier-Stokes in 2D is globally well-posed. Since the current model has some added retardation effect (with the mutual friction), one can argue heuristically that this cannot hurt the extension of global regularity to this nonlinearly coupled two-fluid system (see remarks following the main results). This is the physical intuition behind the main results stated below. \\

\begin{thm} [Global well-posedness for high-regularity data] \label{GWP high regularity}
    Let $m>\frac{d}{2}+1 = 2$, and $u_n^0,u_s^0 \in H^m(\R^2)$ be a pair of divergence-free initial velocity fields. Then, 
    \begin{enumerate}[(i)]
        \item there exist a unique pair of mild solutions $u_n,u_s \in C([0,\infty[ ; H^m(\R^2)) \cap \\ L^2([0,\infty[ ; \dot{H}^{m+1}(\R^2))$ to the incompressible, isothermal 2D HVBK equations described by \eqref{incompressibility}-\eqref{mutual friction}; and \\
        
        \item the solution map is Lipschitz continuous from $H^m(\R^2)\mapsto C([0,\infty[ ; H^m(\R^2)) \cap L^2([0,\infty[ ; \dot{H}^{m+1}(\R^2))$. \\
        
        \item In addition, if the Sobolev index is upgraded to $m\ge 4$, then there also exist a unique pair of pressure fields $p_n,p_s \in C([0,\infty[ ; H^m(\R^2))$; \\
        
        \item The fields are instantaneously smoothed (due to viscosity), so that $u_n,u_s,p_n,p_s \in C^{\infty}(]0,\infty[ ; H^{\infty}(\R^2))$.
    \end{enumerate}
\end{thm}

The proof of Theorem \ref{GWP high regularity} will utilise a contraction mapping argument for local well-posedness. A Beale-Kato-Majda \cite{Beale1984RemarksEquations} analogue will be used to show global well-posedness, and elliptic regularity applied to the pressure-Poisson equation for the pressure field estimates. Once global well-posedness is established for high-regularity data, we can easily extend the result to the following corollary. \\

\begin{cor} [Global well-posedness for smooth, compactly supported data] \label{GWP smooth data}
    Given a pair of divergence-free velocity fields $u_n^0,u_s^0 \in C^{\infty}_c(\R^2)$, there exist a unique pair of smooth (classical) solutions $u_n,u_s \in C^{\infty}([0,\infty[;H^{\infty}(\R^2))$ to the incompressible, isothermal 2D HVBK equations described by \eqref{incompressibility}-\eqref{mutual friction}. Moreover, the pressure fields are also smooth and unique, i.e., $p_n,p_s \in C^{\infty}([0,\infty[;H^{\infty}(\R^2))$.
\end{cor}

\medskip

At this stage, some remarks are in order, particularly about the mutual friction term and the difficulties it poses: 
\begin{enumerate}
    \item Superfluid He-4 is completely inviscid (and dissipation-free) below its critical velocity. As mentioned in Section \ref{intro}, we consider a viscous fluid to help with the analysis. The inviscid-superfluid version of the HVBK equations have not been studied, either numerically or analytically, to the best of the authors' knowledge. \\
    
    \item The physical intuition that the addition of a dissipative mutual friction to the 2D equations would not affect their global regularity need not be valid when the superfluid is treated as inviscid. In fact, it is not known whether the system is even locally well-posed in this case. Note that while the advective nonlinearity preserves the $L^2_x$ norm of the vorticity, the same is not true of the nonlinear mutual friction. This means that the vorticity may not be bounded in $L^{\infty}_t L^2_x$, possibly leading to a breakdown of the system. This is still an open, and interesting, problem. \\
    
    \item The ability of the advective term to be written as the gradient of a symmetric tensor $\left[ \nabla\cdot (u\otimes u) \right]$, also known as the \textit{conservation form}, enables easier manipulation of the nonlinearity. For instance, one of the most useful properties of high-regularity Sobolev spaces $\left( H^m_x, m>\frac{d}{2} \right)$ is that they form an algebra so that $\lVert uv \rVert_{H^m_x} \lesssim \lVert u \rVert_{H^m_x}\lVert v \rVert_{H^m_x}$. This fact when applied to the advective term in the conservation form gives us the same contribution from each velocity in the nonlinearity, as opposed to one being the norm of the velocity and the other a norm of the velocity gradient. This is particularly useful when deriving the $L^{\infty}_t L^{\infty}_x$ controlling norm in the analysis of Navier-Stokes (used to show that the existence time for $C^0_t H^m_x$ solutions is independent of the index $m$). The same goes for heat kernel estimates also, where the conservation form allows for the gradient to be accounted for by the heat kernel instead of the velocity fields (albeit at the cost of temporally singular terms). These steps/methods are no longer as effective when it comes to the mutual friction, since it cannot be reduced to conservation form. There is a clear demarcation, heuristically, of the velocity and the velocity gradient (or vorticity) factors in this term. This asymmetry in the nonlinearity, along with the lopsided presence of only the superfluid vorticity, induces difficulties in the treatment of the mutual friction term. This would be especially true when considering the inviscid superfluid $-$ the lack of a conservation form of mutual friction means one cannot get around the analysis of the gradient of the velocity (more precisely, vorticity), a luxury that is not afforded when lacking viscous diffusion. \\
    
    \item Another problem associated with the mutual friction term is related to the estimates of the pressure. When seeking mild solutions, we eliminate the pressure gradient term from the Navier-Stokes using the Leray projector. This gives us an evolution equation for the velocity field alone. Once a mild solution is established, we can get some estimates on the pressure field by taking the divergence of the Navier-Stokes and using incompressibility. This gives us the \textit{pressure-Poisson equation}.
    
    \begin{equation} \label{pressure-Poisson equation}
        -\Delta p =\nabla\cdot((u\cdot\nabla) u) = \nabla\cdot\nabla\cdot (u\otimes u)
    \end{equation}
    
    This means $p = (-\Delta)^{-1}\nabla\cdot\nabla\cdot (u\otimes u)$. Note how the derivatives roughly ``cancel" out on the RHS. Thus, $p$ is expected to have a regularity similar to $u\otimes u$. This is indeed true (at least in one direction) $-$ since the kernel that maps $u\otimes u$ to $p$ is rotationally symmetric and decays as $|x|^{-d}$, we can utilise the Calderon-Zygmund inequality to conclude $\lVert p \rVert_{L^q_x} \lesssim \lVert u\otimes u \rVert_{L^q_x} \le \lVert u \rVert^2_{L^{2q}_x} \text{ for } 1<q<\infty$. This extension to almost all Lebesgue spaces would not be possible without the conservation form of the advective term. With mutual friction added, there are extra terms on the RHS of \eqref{pressure-Poisson equation}. Clearly, these terms do not have a conservation form; as a result, they do not admit such a simple estimate for the pressure. There is still some respite though: when we have smooth solutions (as in Corollary \ref{GWP smooth data}). In this case, the pressure field is also smooth (see Section \ref{GWP for smooth data section}). Also, as stated in Theorem \ref{GWP high regularity}, if the assumption on the Sobolev index $m$ is upgraded to $m\ge 4$, then we have unique $p_n,p_s \in C^0_t H^m_x$. The reason for this upgrade is again due to the non-conservation form of the mutual friction, which means the RHS in the pressure-Poisson equation contains some terms involving the gradient of the vorticity (equivalently, two derivatives of the velocity).
\end{enumerate}

\section{Apriori estimates} \label{apriori estimates}
In this section, we will derive some apriori estimates for the kinetic energy and enstrophy. 

\subsection{Energy estimate}
Multiplying the normal fluid equation \eqref{normal fluid} by $\rho_n u_n$, the superfluid equation \eqref{superfluid} by $\rho_s u_s$, adding the resulting equations, and integrating over $\R^2$, we obtain (after integrating by parts with vanishing fields at infinity, and using the incompressibility conditions \eqref{incompressibility}):

\begin{align} \label{energy estimate}
    \frac{d}{dt}\left( \frac{1}{2}\rho_n \lVert u_n \rVert_{L^2_x}^2 + \frac{1}{2}\rho_s \lVert u_s \rVert_{L^2_x}^2\right) &+ \rho_n\nu_n \lVert \nabla u_n \rVert_{L^2_x}^2 + \rho_s\nu_s \lVert \nabla u_s \rVert_{L^2_x}^2 \nonumber \\
    &+ \frac{\rho_n \rho_s}{\rho}\frac{B}{2} \int_{\R^2}|\omega_s||u_n-u_s|^2 = 0
\end{align}

Dropping the last (non-negative) term, and integrating over time, we get:

\begin{align} \label{energy and dissipation bound}
    \left( \frac{1}{2}\rho_n \lVert u_n \rVert_{L^\infty_t L^2_x}^2 + \frac{1}{2}\rho_s \lVert u_s \rVert_{L^\infty_t L^2_x}^2\right) &+ \rho_n\nu_n \lVert \nabla u_n \rVert_{L^2_t L^2_x}^2 + \rho_s\nu_s \lVert \nabla u_s \rVert_{L^2_t L^2_x}^2 \nonumber \\
    &\le \frac{1}{2}\rho_n \lVert u_n^0 \rVert_{L^2_x}^2 + \frac{1}{2}\rho_s \lVert u_s^0 \rVert_{L^2_x}^2
\end{align}

Thus, we see that the energy at any time (and dissipation \textit{upto} any time) is bounded above by the initial energy. Recall that we will choose initial data satisfying the assumptions of Lemma \ref{finite kinetic energy}, so that the kinetic energy is finite at $t=0$.

\subsection{Vorticity equations} \label{vorticity equations}
Operating on \eqref{normal fluid} and \eqref{superfluid} with the curl operator, we get the evolution equations for the scalar vorticity fields (which are also incompressible).

\begin{gather}
    \partial_t \omega_n + (u_n\cdot \nabla)\omega_n = \nu_n \Delta \omega_n + \frac{\rho_s}{\rho}\T \label{normal fluid vorticity} \\
    \partial_t \omega_s + (u_s\cdot \nabla)\omega_s = \nu_s \Delta \omega_s - \frac{\rho_n}{\rho}\T \label{superfluid vorticity} \\
    \mathcal{T} = -\frac{B}{2}\nabla\times [|\omega_s|(u_n - u_s)]  + \frac{B'}{2}(u_n - u_s)\cdot\nabla\omega_s \label{mutual friction torque}
\end{gather}

\subsection{Enstrophy estimate} \label{enstrophy estimate section}
Just as in section \ref{energy estimate}, we multiply \eqref{normal fluid vorticity} by $\rho_n \omega_n$, \eqref{superfluid vorticity} by $\rho_s \omega_s$, add the resulting equations, and integrate over $\R^2$ to arrive at:

\begin{align} \label{enstrophy estimate}
    \frac{d}{dt}\left( \frac{1}{2}\rho_n \lVert \omega_n \rVert_{L^2_x}^2 + \frac{1}{2}\rho_s \lVert \omega_s \rVert_{L^2_x}^2 \right) &+ \rho_n\nu_n \lVert \nabla \omega_n \rVert_{L^2_x}^2 + \rho_s\nu_s \lVert \nabla \omega_s \rVert_{L^2_x}^2 \nonumber \\
    &+ \frac{\rho_n \rho_s}{\rho}\frac{B}{2} \int_{\R^2}|\omega_s||\omega_n - \omega_s|^2 \nonumber \\
    &= \frac{\rho_n \rho_s}{2\rho} \int_{\R^2} (\omega_n - \omega_s) (u_n - u_s)\cdot \left( B'\nabla\omega_s - B\nabla^{\perp}\lvert\omega_s\rvert \right)
\end{align}

Dropping the last term (non-negative) on the LHS, and using H\"older's inequality,

\begin{align*}
    \frac{d}{dt}\left( \frac{1}{2}\rho_n \lVert \omega_n \rVert_{L^2_x}^2 + \frac{1}{2}\rho_s \lVert \omega_s \rVert_{L^2_x}^2 \right) &+ \rho_n\nu_n \lVert \nabla \omega_n \rVert_{L^2_x}^2 + \rho_s\nu_s \lVert \nabla \omega_s \rVert_{L^2_x}^2 \\
    &\le \frac{\rho_n \rho_s}{2\rho}(B+B') \lVert \omega_n - \omega_s \rVert_{L^4_x} \lVert u_n - u_s \rVert_{L^4_x} \lVert \nabla\omega_s \rVert_{L^2_x}
\end{align*}

We use Cauchy's inequality $\left( ab\le \epsilon a^2 + \frac{b^2}{4\epsilon} \right)$ to extract the $\lVert \nabla\omega_s \rVert_{L^2_x}$ term from the RHS and absorb into the LHS. We also use the two-dimensional Ladyzhenskaya inequality $\left( \lVert u \rVert^2_{L^4_x} \lesssim \lVert u \rVert_{L^2_x}\lVert \nabla u \rVert_{L^2_x} \right)$ fo handling the other terms on the RHS. This leaves us with:

\begin{align*}
    \frac{d}{dt}\left( \frac{1}{2}\rho_n \lVert \omega_n \rVert_{L^2_x}^2 + \frac{1}{2}\rho_s \lVert \omega_s \rVert_{L^2_x}^2 \right) &+ \rho_n\nu_n \lVert \nabla \omega_n \rVert_{L^2_x}^2 + \frac{\rho_s\nu_s}{2} \lVert \nabla \omega_s \rVert_{L^2_x}^2 \\
    &\lesssim \lVert u_n - u_s \rVert_{L^2_x} \lVert \nabla u_n - \nabla u_s \rVert_{L^2_x} \lVert \omega_n - \omega_s \rVert_{L^2_x} \lVert \nabla\omega_n - \nabla\omega_s \rVert_{L^2_x}
\end{align*}

Using the $L^{\infty}_t L^2_x$ bound on $u_n$ and $u_s$ from \eqref{energy and dissipation bound}, and Cauchy's inequality again, the RHS can be controlled as follows.

\begin{align*}
    RHS &\lesssim \lVert \omega_n - \omega_s \rVert^2_{L^2_x} \lVert \nabla u_n - \nabla u_s \rVert^2_{L^2_x} + \frac{\rho_n\nu_n}{2}\lVert \nabla\omega_n \rVert^2_{L^2_x} + \frac{\rho_s\nu_n}{4}\lVert \nabla\omega_s \rVert^2_{L^2_x}
\end{align*}

The second and third terms on the RHS can be absorbed into the corresponding terms on the LHS. 

\begin{equation} \label{enstrophy balance equation}
    \begin{split}
        \frac{d}{dt}\left( \frac{\rho_n}{2} \lVert \omega_n \rVert_{L^2_x}^2 + \frac{\rho_s}{2} \lVert \omega_s \rVert_{L^2_x}^2 \right) &+ \frac{\rho_n\nu_n}{2} \lVert \nabla \omega_n \rVert_{L^2_x}^2 + \frac{\rho_s\nu_s}{4} \lVert \nabla \omega_s \rVert_{L^2_x}^2 \\
        &\lesssim \lVert \omega_n - \omega_s \rVert^2_{L^2_x} \lVert \nabla u_n - \nabla u_s \rVert^2_{L^2_x} \\
        &\lesssim \left( \frac{\rho_n}{2} \lVert \omega_n \rVert_{L^2_x}^2 + \frac{\rho_s}{2} \lVert \omega_s \rVert_{L^2_x}^2 \right) \left[\lVert \nabla u_n \rVert^2_{L^2_x} + \lVert \nabla u_s \rVert^2_{L^2_x}\right]
    \end{split}
\end{equation}

We then drop non-negative terms on the LHS, and use Gr\"onwall's inequality to arrive at:

\begin{align}
    \frac{\rho_n}{2} \lVert \omega_n \rVert_{L^2_x}^2(t) + \frac{\rho_s}{2} \lVert \omega_s \rVert_{L^2_x}^2(t) &\le \left( \frac{\rho_n}{2} \lVert \omega^0_n \rVert_{L^2_x}^2 + \frac{\rho_s}{2} \lVert \omega^0_s \rVert_{L^2_x}^2 \right) e^{C\left[\lVert \nabla u_n \rVert^2_{L^2_t L^2_x} + \lVert \nabla u_s \rVert^2_{L^2_t L^2_x}\right]}
\end{align}

Since $\lVert \nabla u \rVert_{L^2_t L^2_x}$ is bounded for each fluid according to \eqref{energy and dissipation bound}, we observe that $\omega_n$ and $\omega_s$ are bounded in $L^{\infty}_t L^2_x$. Substituting this upper bound on the vorticities back into \eqref{enstrophy balance equation}, we integrate over time to calculate upper bounds for $\lVert \nabla \omega_n \rVert_{L^2_t L^2_x}$ and $\lVert \nabla \omega_s \rVert_{L^2_t L^2_x}$ as well. In summary,

\begin{equation} \label{enstrophy and palinstrophy bounds}
    \lVert \omega_n \rVert_{L^{\infty}_t L^2_x},\lVert \omega_s \rVert_{L^{\infty}_t L^2_x},\lVert \nabla \omega_n \rVert_{L^2_t L^2_x},\lVert \nabla \omega_s \rVert_{L^2_t L^2_x} \le C\left( \lVert \omega^0_n \rVert_{L^2_x},\lVert \omega^0_s \rVert_{L^2_x},\lVert u^0_n \rVert_{L^2_x},\lVert u^0_s \rVert_{L^2_x} \right)
\end{equation}

The apriori bounds in \eqref{energy and dissipation bound} and \eqref{enstrophy and palinstrophy bounds} will be repeatedly used in the proof of global well-posedness (section \ref{global well-posedness}).

\section{Local well-posedness} \label{local well-posedness}
In this section, we will establish the local well-posedness of the system \eqref{incompressibility}-\eqref{mutual friction} for high-regularity data. We will use the properties of the heat kernel to set up a contraction mapping that allows us to show the existence of a unique (mild) solution. Since we have a two-fluid system, the contraction map will be ``alternating", in the sense that for a given $u_s$, we contract the normal fluid map, and for a given $u_n$, we contract the superfluid map. \\

The appropriate space for the contraction is denoted by $X_i$, and the corresponding norm is defined as:

\begin{equation} \label{The Norm}
    \lVert u_i \rVert_{X_i} := \lVert u_i \rVert_{C^0_t H^m_x} + \nu_i^{\frac{1}{2}} \lVert \nabla u_i \rVert_{L^2_t H^m_x} 
\end{equation}

\noindent where $i\in \{n,s\}$. The space-time variables have the domains $t\in [0,T[$ and $x\in\R^2$.

\subsection{Mild solution and bound on the norm} \label{Mild solution and bound on the norm}
Starting with divergence-free initial data (in accordance with Lemma \ref{finite kinetic energy}) $u_n^0,u_s^0\in H^m(\R^2)$ for $m>\frac{d}{2}+1=2$, we write the required contraction maps in the form of mild solutions corresponding to \eqref{normal fluid} and \eqref{superfluid}.

\begin{gather}
    \Phi(t) = e^{\nu_n t\Delta}u_n^0 - \int_0^t e^{\nu_n (t-\tau)\Delta}\Leray\left[u_n\cdot\nabla u_n + \beta_n |\omega_s|(u_n - u_s) - \beta_n' \omega_s \times (u_n - u_s) \right] \label{Phi map} \\
    \Psi(t) = e^{\nu_s t\Delta}u_s^0 - \int_0^t e^{\nu_s (t-\tau)\Delta}\Leray\left[u_s\cdot\nabla u_s - \beta_s |\omega_s|(u_n - u_s) + \beta_s' \omega_s \times (u_n - u_s) \right] \label{Psi map}
\end{gather}

Here, $u_n,u_s\in H^m(\R^2)$ as well, and the (positive) constants $\beta$'s are defined below.

\begin{gather*}
    \beta_n = \frac{\rho_s}{\rho}\frac{B}{2} \quad , \quad \beta_n' = \frac{\rho_s}{\rho}\frac{B'}{2} \\
    \beta_s = \frac{\rho_n}{\rho}\frac{B}{2} \quad , \quad \beta_s' = \frac{\rho_n}{\rho}\frac{B'}{2}
\end{gather*}

From \eqref{Phi map}, we have almost everywhere (a.e.) in time:

\begin{equation*}
    \frac{\partial \Phi}{\partial t} = \nu_n \Delta\Phi + F_n
\end{equation*}

where $F_n = -\Leray\left[u_n\cdot\nabla u_n + \beta_n |\omega_s|(u_n - u_s) - \beta_n' \omega_s \times (u_n - u_s) \right]$. \\

For each multi-index $\alpha$ with $0\le|\alpha|\le m$, we can derive a higher-order energy inequality for the contraction maps in \eqref{Phi map} and \eqref{Psi map}.

\begin{equation*}
    \frac{\partial}{\partial t} D^{\alpha}\Phi = \nu_n \Delta D^{\alpha}\Phi + D^{\alpha}F_n 
\end{equation*}

Multiplying the above equation by $D^{\alpha}\Phi$, integrating over $\R^2$, using the Holder's inequality, and summing over all $0\le |\alpha|\le m$ leads to:

\begin{equation*}
    \frac{d}{dt}\lVert \Phi \rVert^2_{H^m_x} + \nu_n \lVert \nabla\Phi \rVert^2_{H^m_x} \le \lVert \Phi \rVert_{H^m_x} \lVert F_n \rVert_{H^m_x}
\end{equation*}

Integrating over time, from $0$ to $T$, and using the fact that $\Phi(0) = u_n^0$:

\begin{equation*}
    \lVert \Phi \rVert^2_{C^0_t H^m_x} + \nu_n \lVert \nabla\Phi \rVert^2_{L^2_t H^m_x} \le \lVert u_n^0 \rVert^2_{H^m_x} + \lVert \Phi \rVert_{C^0_t H^m_x} \lVert F_n \rVert_{L^1_t H^m_x}
\end{equation*}

Finally, we use the inequality $ab\le \frac{a^2 + b^2}{2}$ to simplify the second term on the RHS and obtain:

\begin{gather}
    \lVert \Phi \rVert_{X_n} \lesssim \lVert u_n^0 \rVert_{H^m_x} + \lVert F_n \rVert_{L^1_t H^m_x} \label{Phi map norm} \\
    \lVert \Psi \rVert_{X_s} \lesssim \lVert u_s^0 \rVert_{H^m_x} + \lVert F_s \rVert_{L^1_t H^m_x} \label{Psi map norm}
\end{gather}

Equation \eqref{Psi map norm} is derived in the same way as \eqref{Phi map norm}.

\subsection{The contraction} \label{The contraction}

Since Sobolev spaces with $m>\frac{d}{2}$ form an algebra, i.e., $\lVert uv \rVert_{H^m_x} \lesssim \lVert u \rVert_{H^m_x} \lVert v \rVert_{H^m_x}$,

\begin{equation*}
    \lVert F_n \rVert_{H^m_x} \lesssim \lVert u_n \rVert_{H^m_x} \lVert \nabla u_n \rVert_{H^m_x} + (\beta_n + \beta_n') \lVert \omega_s \rVert_{H^m_x}\lVert u_n - u_s \rVert_{H^m_x}
\end{equation*}

From the Calderon-Zygmund inequality, we have $\lVert \omega_s \rVert_{H^m_x} \lesssim \lVert \nabla u_s \rVert_{H^m_x}$.  Integrating over the time interval $[0,T_n[$, 

\begin{align}
    \lVert F_n \rVert_{L^1_t H^m_x} &\lesssim \lVert u_n \rVert_{C^0_t H^m_x} \lVert \nabla u_n \rVert_{L^2_t H^m_x} T_n^{\frac{1}{2}} + (\beta_n + \beta_n') \lVert \nabla u_s \rVert_{L^2_t H^m_x}\lVert u_n - u_s \rVert_{C^0_t H^m_x} T_n^{\frac{1}{2}} \nonumber \\
    &\lesssim \nu_n^{-\frac{1}{2}} T_n^{\frac{1}{2}} \lVert u_n \rVert^2_{X_n} + (\beta_n + \beta_n') \nu_s^{-\frac{1}{2}} T_n^{\frac{1}{2}} \lVert u_s \rVert_{X_s} \left( \lVert u_n \rVert_{X_n} + \lVert u_s \rVert_{X_s} \right) \label{F_n norm}
\end{align}

From \eqref{F_n norm} and \eqref{Phi map norm}, and repeating the same procedure for $\Psi$, we obtain:

\begin{gather*}
    \lVert \Phi \rVert_{X_n} \lesssim \lVert u_n^0 \rVert_{H^m_x} + T_n^{\frac{1}{2}} \left[ \nu_n^{-\frac{1}{2}} \lVert u_n \rVert^2_{X_n} + \nu_s^{-\frac{1}{2}} (\beta_n + \beta_n') \lVert u_s \rVert_{X_s} \left( \lVert u_n \rVert_{X_n} + \lVert u_s \rVert_{X_s} \right) \right] \\
    \lVert \Psi \rVert_{X_s} \lesssim \lVert u_s^0 \rVert_{H^m_x} + T_s^{\frac{1}{2}} \left[ \nu_s^{-\frac{1}{2}} \lVert u_s \rVert^2_{X_s} + \nu_s^{-\frac{1}{2}} (\beta_s + \beta_s') \lVert u_s \rVert_{X_s} \left( \lVert u_n \rVert_{X_n} + \lVert u_s \rVert_{X_s} \right) \right]
\end{gather*}

Let $B_X (R)$ denote a ball (centred at the origin) of radius $R$ in the space $X$. Consider $u_n \in B_{X_n}(R_n)$ and $u_s \in B_{X_s}(R_s)$. Thus, the above equations therefore simplify to:

\begin{gather}
    \lVert \Phi \rVert_{X_n} \lesssim \lVert u_n^0 \rVert_{H^m_x} + T_n^{\frac{1}{2}} \left[ \nu_n^{-\frac{1}{2}} R_n^2 + \nu_s^{-\frac{1}{2}} (\beta_n + \beta_n') R_s \left( R_n + R_s \right) \right] \label{Phi map norm updated} \\
    \lVert \Psi \rVert_{X_s} \lesssim \lVert u_s^0 \rVert_{H^m_x} + T_s^{\frac{1}{2}} \left[ \nu_s^{-\frac{1}{2}} R_s^2 + \nu_s^{-\frac{1}{2}} (\beta_s + \beta_s') R_s \left( R_n + R_s \right) \right] \label{Psi map norm updated}
\end{gather}

Now, we will fix the superfluid velocity field $u_s \in B_{X_s}(R_s)$ and the initial normal fluid velocity $u_n^0$, and consider two different normal fluid velocity fields $u_n^1, u_n^2 \in B_{X_n}(R_n)$. Then the difference of the $\Phi$ map for the two different velocities is:

\begin{equation*}
    \Phi(u_n^1) - \Phi(u_n^2) = \int_0^t e^{\nu_n (t-\tau)\Delta} \left( F_n(u_n^1) - F_n(u_n^2) \right) d\tau
\end{equation*}

Using the estimate from \eqref{Phi map norm},

\begin{align*}
    \Rightarrow \lVert \Phi(u_n^1) - \Phi(u_n^2) \rVert_{X_n} &\lesssim \lVert F_n^1 - F_n^2 \rVert_{L^1_t H^m_x} \\
    &\lesssim \lVert u_n^1\cdot\nabla u_n^1 - u_n^2\cdot\nabla u_n^2 - \beta_n |\omega_s|(u_n^1 - u_n^2) + \beta_n' \omega_s \times (u_n^1 - u_n^2) \rVert_{L^1_t H^m_x} \\
    &\lesssim T_n^{\frac{1}{2}} \left[ \nu_n^{-\frac{1}{2}}(\lVert u_n^1 \rVert_{X_n} + \lVert u_n^2 \rVert_{X_n}) + (\beta_n + \beta_n') \nu_s^{-\frac{1}{2}} \lVert u_s \rVert_{X_s} \right] \lVert u_n^1 - u_n^2 \rVert_{X_n}
\end{align*}

A similar procedure can be repeated for the $\Psi$ map, starting with the initial data $u_s^0$, fixing the normal fluid velocity field $u_n \in B_{X_n}(R_n)$, and considering two different superfluid velocity\footnote{Note that in this case, $\omega_s$ is dependent on which of the superfluid velocity fields $u_s^1$ or $u_s^2$ is in consideration in that term.} fields $u_s^1, u_s^2 \in B_{X_s}(R_s)$. Thus,

\begin{gather}
    \lVert \Phi(u_n^1) - \Phi(u_n^2) \rVert_{X_n} \lesssim T_n^{\frac{1}{2}} \left[ \nu_n^{-\frac{1}{2}}R_n + (\beta_n + \beta_n') \nu_s^{-\frac{1}{2}} R_s \right] \lVert u_n^1 - u_n^2 \rVert_{X_n} \label{Phi Lipschitz contraction} \\
    \lVert \Psi(u_s^1) - \Psi(u_s^2) \rVert_{X_s} \lesssim T_s^{\frac{1}{2}} \left[ \nu_s^{-\frac{1}{2}}(1 + \beta_s + \beta_s')R_s + (\beta_s + \beta_s') \nu_s^{-\frac{1}{2}} R_n \right] \lVert u_s^1 - u_s^2 \rVert_{X_n} \label{Psi Lipschitz contraction}
\end{gather}

From \eqref{Phi map norm updated}, \eqref{Psi map norm updated}, \eqref{Phi Lipschitz contraction} and \eqref{Psi Lipschitz contraction}, we see that if the following (sufficient) conditions are satisfied, the maps $\Phi$ and $\Psi$ will be contractions. ($c_1,c_2,c_3$ and $c_4$ are positive constants that were suppressed in all the inequalities so far.)

\begin{align}
    c_1 \left[ \lVert u_n^0 \rVert_{H^m_x} + T_n^{\frac{1}{2}} \left( \nu_n^{-\frac{1}{2}} R_n^2 + \nu_s^{-\frac{1}{2}} (\beta_n + \beta_n') R_s \left( R_n + R_s \right) \right) \right] &\le R_n \label{Phi map into} \\
    c_2 \left[ T_n^{\frac{1}{2}} \left( \nu_n^{-\frac{1}{2}}R_n + (\beta_n + \beta_n') \nu_s^{-\frac{1}{2}} R_s \right) \right] &< 1 \label{Phi contraction} \\
    c_3 \left[ \lVert u_s^0 \rVert_{H^m_x} + T_s^{\frac{1}{2}} \left( \nu_s^{-\frac{1}{2}} R_s^2 + \nu_s^{-\frac{1}{2}} (\beta_s + \beta_s') R_s \left( R_n + R_s \right) \right) \right] &\le R_s \label{Psi map into} \\
    c_4 \left[ T_s^{\frac{1}{2}} \left( \nu_s^{-\frac{1}{2}}(1 + \beta_s + \beta_s')R_s + (\beta_s + \beta_s') \nu_s^{-\frac{1}{2}} R_n \right) \right] &< 1 \label{Psi contraction}
\end{align}

Define $\delta = 1+\beta_n + \beta_n' + \beta_s + \beta_s'$ and $\nu^{-\frac{1}{2}} = \nu_n^{-\frac{1}{2}} + \nu_s^{-\frac{1}{2}}$. Then, the conditions in \eqref{Phi map into}-\eqref{Psi contraction} are automatically satisfied if the following more conservative inequalities are obeyed (for some sufficiently large $N_n,N_s\in \N$).

\begin{align}
    c_1 \delta \left[ \lVert u_n^0 \rVert_{H^m_x} + T_n^{\frac{1}{2}} \nu^{-\frac{1}{2}}(R_n + R_s)^2 \right] &\le R_n \label{Phi map into conservative} \\
    c_2 \delta T_n^{\frac{1}{2}} \nu^{-\frac{1}{2}}(R_n + R_s) &\le \frac{1}{N_n} \label{Phi contraction conservative} \\
    c_3 \delta \left[ \lVert u_s^0 \rVert_{H^m_x} + T_s^{\frac{1}{2}} \nu^{-\frac{1}{2}}(R_n + R_s)^2 \right] &\le R_s \label{Psi map into conservative} \\
    c_4 \delta T_s^{\frac{1}{2}} \nu^{-\frac{1}{2}}(R_n + R_s) &\le \frac{1}{N_s} \label{Psi contraction conservative}
\end{align}

Upper bounds for $T_n$ and $T_s$ from \eqref{Phi contraction conservative} and \eqref{Psi contraction conservative} are substituted in \eqref{Phi map into conservative} and \eqref{Psi map into conservative} to give:

\begin{equation*}
    c_1 \delta \lVert u_n^0 \rVert_{H^m_x} + c_3 \delta \lVert u_s^0 \rVert_{H^m_x} + \left( \frac{c_1}{N_n c_2} + \frac{c_3}{N_s c_4} \right) (R_n + R_s) \le (R_n + R_s)
\end{equation*}

This is easily satisfied by choosing $R_n = 2c_1 \delta \lVert u_n^0 \rVert_{H^m_x}$ and $R_s = 2c_3 \delta \lVert u_s^0 \rVert_{H^m_x}$, and sufficiently large $N_n$ and $N_s$. This choice also leads to the conclusion that:

\begin{gather}
    T_n \le \frac{\nu}{(2c_2 \delta^2N_n)^2}\frac{1}{(c_1\lVert u_n^0 \rVert_{H^m_x} + c_3\lVert u_s^0 \rVert_{H^m_x})^2} \label{T_n upper bound} \\
    T_s \le \frac{\nu}{(2c_2 \delta^2N_s)^2}\frac{1}{(c_1\lVert u_n^0 \rVert_{H^m_x} + c_3\lVert u_s^0 \rVert_{H^m_x})^2} \label{T_s upper bound}
\end{gather}

From \eqref{Phi contraction conservative}-\eqref{Psi contraction conservative}, we observe that the upper bounds on the existence times are inversely proportional to $R_n+R_s$. This allows us to use a bootstrapping argument to establish that there exists a \textit{unique maximal time} $T_*\in(0,\infty]$ such that $\lVert u_n \rVert_{X_n}(t) \lVert u_n \rVert_{X_s}(t) \rightarrow \infty$ as $t\rightarrow T_*$. In other words, at least one of the two velocity fields blows up in the respective $X_i$ norm as the maximal time is approached.

\subsection{Lipschitz continuous dependence on initial data}
So far, we have shown the local existence and uniqueness of (mild) solutions. We will now establish Lipschitz continuous dependence of the solutions on initial data, to complete the well-posedness proof. Consider two initial data $u_n^0, v_n^0$ with corresponding mild solutions $u_n,v_n$ (for a fixed $u_s$). Then,

\begin{align*}
    \lVert u_n - v_n\rVert_{X_n} \lesssim \lVert u_n^0 - v_n^0\rVert_{H^m_x} + \lVert u_n\cdot\nabla u_n - v_n\cdot\nabla v_n &- \beta_n |\omega_s|(u_n - v_n) \\ 
    &+ \beta_n' \omega_s\times(u_n - v_n) \rVert_{L^1_t H^m_x}
\end{align*}

Just as done in section \ref{The contraction}, we can show:

\begin{multline*}
    \lVert u_n - v_n\rVert_{X_n} \lesssim \lVert u_n^0 - v_n^0\rVert_{H^m_x} \\ 
    + T_*^{\frac{1}{2}}\left[ \nu_n^{-\frac{1}{2}}(\lVert u_n \rVert_{X_n} + \lVert v_n \rVert_{X_n}) + \nu_s^{-\frac{1}{2}}\lVert u_s \rVert_{X_s} \right]\lVert u_n - v_n \rVert_{X_n}
\end{multline*}

From \eqref{Phi Lipschitz contraction}, \eqref{Phi contraction} and \eqref{Phi contraction conservative}, we can choose $T_*$ sufficiently small so that the second term on the RHS may be absorbed into the LHS. This gives the required result for the normal fluid. The same procedure can be repeated for the superfluid as well.

\subsection{Regularity of the pressure fields} \label{regularity of pressure}

We will now briefly comment on the regularity of the pressure fields. Applying the divergence operator on \eqref{normal fluid} and \eqref{superfluid}, and using incompressibility, we get the Poisson equations that govern the evolution of the pressure fields.

\begin{equation} \label{pressure Poisson equations}
    \begin{split}
        -\Delta p_n &= \nabla\cdot\nabla\cdot(u_n \otimes u_n) + \beta_n (u_n - u_s)\cdot\nabla|\omega_s| - \beta_n' \nabla\cdot (\omega_s \times (u_n - u_s)) \\
        -\Delta p_s &= \nabla\cdot\nabla\cdot(u_s \otimes u_s) - \beta_s (u_n - u_s)\cdot\nabla|\omega_s| + \beta_s' \nabla\cdot (\omega_s \times (u_n - u_s))
    \end{split}
\end{equation}

\medskip

The solutions $p_n,p_s$ of these equations must be unique. At first look, it appears that an arbitrary harmonic function could be added to the pressure fields, breaking uniqueness. But since we seek bounded pressures, Liouville's theorem (see chapter 2 of \cite{EvansLawrenceC.2010PartialEquations}) guarantees that such a harmonic function is in fact a constant. This constant can be forced to be zero, by demanding that the pressure vanishes at infinity. \\

Now, for unique solutions of the Poisson equation $-\Delta p = f$, with $p =0$ on the boundary, elliptic regularity theory estimates (see chapter 6 of \cite{EvansLawrenceC.2010PartialEquations}) show that $\lVert p \rVert_{H^{m+2}_x} \lesssim \lVert f \rVert_{H^m_x}$. The RHS of each equation in \eqref{pressure Poisson equations} consists of terms like $(\nabla u)(\nabla u)$, and $(u) (\nabla \omega) \sim (u) (D^2 u)$. If we upgrade our Sobolev index to $m\ge 4$, then we have $m-2 > \frac{d}{2}$, allowing us to use the algebra property of $H^{m-2}_x$. Thus, 

\begin{align*}
    \lVert \text{RHS} \rVert_{H^{m-2}_x} &\lesssim \lVert (\nabla u)(\nabla u) \rVert_{H^{m-2}_x} + \lVert (u)(D^2 u) \rVert_{H^{m-2}_x} \\
    &\lesssim \lVert (\nabla u) \rVert^2_{H^{m-2}_x} + \lVert u \rVert_{H^{m-2}_x}\lVert D^2 u \rVert_{H^{m-2}_x} \\
    &\lesssim \lVert u \rVert^2_{H^{m}_x} < \infty
\end{align*}

Therefore, this shows that the RHS of \eqref{pressure Poisson equations} $\in C^0_t H^{m-2}_x\Rightarrow p_n,p_s\in C^0_t H^m_x$.

\subsection{Instantaneous smoothing} \label{instantaneous smoothing}

Viscosity-driven momentum diffusion leads to the well-known instantaneous smoothing of the solution. This standard result can be argued as follows, where $u$ can mean either the superfluid or the normal fluid velocity. Since we have established that $u \in C^0([0,T_*[;\dot{H}^m_x)\cap L^2([0,T_*[;\dot{H}^{m+1}_x)$, we can conclude that $\lVert \nabla u \rVert_{H^m_x}(t) < \infty$ for almost every $t\in [0,T_*[$. In other words, $\lVert u \rVert_{H^{m+1}_x}(t) < \infty$ for almost every $t\in [0,T_*[$. \\ 

So, for an arbitrary $0< \delta < T_*$, we can find a time $0< t_0 < \delta$ such that $\lVert u \rVert_{H^{m+1}_x}(t_0) < \infty$. Denoting this velocity field as $u^{t_0}$ and considering it as an $H^{m+1}_x$ datum, we can evolve the system (as shown thus far) to find a local $C^0_t H^{m+1}_x$ solution in some time interval $[t_0,T_*'[$, where $T_*' - t_0$ is the maximal existence time corresponding to this $H^{m+1}_x$ datum. The choice of $t_0$ can be made so that $\delta < T_*'$. By the uniqueness of solutions on the interval $[t_0,\delta]$, we may conclude that indeed $u\in C^0 ([t_0,\delta];H^{m+1}_x)$. Since $\delta$ is arbitrary, this implies $u\in C^0 (]0,T_*[;H^{m+1}_x)$. Iterating this argument shows that $u\in C^0 (]0,T_*[;H^{\infty}_x)$. \\

The mild solutions given by \eqref{mild solution normal fluid} and \eqref{mild solution superfluid} satisfy a.e. in time:

\begin{align*}
    \partial_t u_n &= \nu_n \Delta u_n - \Leray (u_n\cdot\nabla u_n + \beta_n |\omega_s|(u_n - u_s) - \beta_n' \omega_s \times (u_n - u_s)) \\
    \partial_t u_s &= \nu_s \Delta u_s - \Leray (u_s\cdot\nabla u_s - \beta_s |\omega_s|(u_n - u_s) + \beta_s' \omega_s \times (u_n - u_s))
\end{align*}

\medskip

Since $u_n,u_s \in C^0(]0,\infty[;H^{\infty}(\R^2))$, each term on the RHS of these equations is in $C^0(]0,\infty[;H^{\infty}(\R^2))$. This means the LHS (the time-derivative term) is in $C^0(]0,\infty[;H^{\infty}(\R^2))$, i.e., $u_n,u_s \in C^1(]0,\infty[;H^{\infty}(\R^2))$. This argument can be iterated to show that indeed $u_n,u_s \in C^{\infty}(]0,\infty[;H^{\infty}(\R^2))$. In other words, the solutions are smooth in time (for $t>0$), and $H^{\infty}$-smooth in space. Repeating the steps of Section \ref{regularity of pressure} proves the uniqueness and space-time smoothness of pressure fields as well (for $t>0$). \\

\qed








\section{Global well-posedness} \label{global well-posedness}

\subsection{Higher-order energy estimate} \label{higher-order energy estimate subsection}

We will now establish global well-posedness for high-regularity data using the method of Beale-Kato-Majda \cite{Beale1984RemarksEquations}. We begin by deriving a higher-order energy estimate like in section \ref{Mild solution and bound on the norm}. Acting on \eqref{normal fluid} with a higher-order derivative operator $D^{\alpha}$ (for $0\le |\alpha| \le m$, multiplying by $D^{\alpha}u_n$, and summing over all $\alpha$, we get:

\begin{align} 
    \frac{d}{dt}\frac{1}{2}\lVert u_n \rVert^2_{H^m_x} + \nu_n\lVert\nabla u_n\rVert^2_{H^m_x} = -\sum_{\alpha}\langle &D^{\alpha}u_n,\Leray D^{\alpha}\left[ u_n\cdot\nabla u_n \right. \nonumber \\ 
    &\left. + \beta_n |\omega_s|(u_n - u_s) - \beta_n' \omega_s\times(u_n - u_s) \right] \rangle \label{normal higher-order preliminary}
\end{align}

\noindent where $\langle , \rangle$ denotes the $L^2_x$ inner product. The Leray projector is self-adjoint and commutes with derivatives, so it can be dropped. By incompressibility, $\langle D^{\alpha}u_n, u_n\cdot\nabla D^{\alpha} u_n \rangle = 0$. So, we add this vanishing quantity to the RHS.

\begin{align*}
    \begin{split}
        RHS &= -\sum_{\alpha}\langle D^{\alpha}u_n,D^{\alpha}\left[ u_n\cdot\nabla u_n \right] - u_n\cdot\nabla D^{\alpha} u_n \rangle \\
        &\qquad \qquad \qquad - \sum_{\alpha}\langle D^{\alpha}u_n ,D^{\alpha} \left[ \beta_n |\omega_s|(u_n - u_s)  - \beta_n' \omega_s\times(u_n - u_s)\right] \rangle
    \end{split}
    \\
    &\lesssim \lVert u_n \rVert_{H^m_x}\left[ \lVert \nabla u_n \rVert_{L^{\infty}_x} \lVert u_n \rVert_{H^m_x} + \lVert |\omega_s|(u_n - u_s) \rVert_{H^m_x} \right] \\
    &\lesssim \lVert u_n \rVert_{H^m_x}\left[ \lVert \nabla u_n \rVert_{L^{\infty}_x} \lVert u_n \rVert_{H^m_x} + \lVert u_n - u_s \rVert_{L^{\infty}_x} \lVert \omega_s \rVert_{H^m_x} + \lVert \omega_s \rVert_{L^{\infty}_x} \lVert u_n - u_s \rVert_{H^m_x} \right]
\end{align*}

\noindent where we used H\"older's inequality and calculus inequalities for Sobolev spaces (Eqns. (3.31) and (3.32) from \cite{Majda2002VorticityFlow}). \\

Using the above estimate of the RHS in \eqref{normal higher-order preliminary},

\begin{align}
    \frac{d}{dt}\frac{1}{2}\lVert u_n \rVert^2_{H^m_x} &+ \nu_n\lVert\nabla u_n\rVert^2_{H^m_x} \nonumber \label{normal higher-order estimate} \\
        & \begin{split} \lesssim \lVert \nabla u_n \rVert_{L^{\infty}_x} \lVert u_n \rVert^2_{H^m_x} &+ \lVert u_n - u_s \rVert_{L^{\infty}_x} \lVert u_n \rVert_{H^m_x} \lVert \omega_s \rVert_{H^m_x} \\
        &+ \lVert \omega_s \rVert_{L^{\infty}_x} \lVert u_n - u_s \rVert_{H^m_x}\lVert u_n \rVert_{H^m_x} \end{split}
\end{align}

A similar calculation for the superfluid equation \eqref{superfluid} can be performed.

\begin{align}
    \frac{d}{dt}\frac{1}{2}\lVert u_s \rVert^2_{H^m_x} &+ \nu_s\lVert\nabla u_s\rVert^2_{H^m_x} \nonumber \label{super higher-order estimate} \\
    & \begin{split} \lesssim \lVert \nabla u_s \rVert_{L^{\infty}_x} \lVert u_s \rVert^2_{H^m_x} &+ \lVert u_n - u_s \rVert_{L^{\infty}_x} \lVert u_s \rVert_{H^m_x} \lVert \omega_s \rVert_{H^m_x} \\
    &+ \lVert \omega_s \rVert_{L^{\infty}_x} \lVert u_n - u_s \rVert_{H^m_x}\lVert u_s \rVert_{H^m_x} \end{split}
\end{align}

Once again, using Cauchy's inequality, we can extract out a $\frac{\nu_s}{2} \lVert \omega_s \rVert^2_{H^m_x}$ from the second terms on the RHS of \eqref{normal higher-order estimate} and \eqref{super higher-order estimate}. Since $\lVert \omega_s \rVert_{H^m_x} \le \lVert \nabla u_s \rVert_{H^m_x}$, we may add the equations and these extracted terms can be used to cancel out the $\nu_s\lVert \nabla u_s \rVert^2_{H^m_x}$ on the LHS. Finally, after dropping $\nu_n \lVert \nabla u_n \rVert^2_{H^m_x}$ on the LHS, the following inequality results. 

\begin{align}
    \frac{d}{dt}\left( \lVert u_n \rVert^2_{H^m_x} + \lVert u_s \rVert^2_{H^m_x} \right) \lesssim &\left( \lVert u_n \rVert^2_{H^m_x} + \lVert u_s \rVert^2_{H^m_x} \right) \times \nonumber \label{higher-order estimate} \\
    &\left[ \lVert \nabla u_n \rVert_{L^{\infty}_x} + \lVert \nabla u_s \rVert_{L^{\infty}_x} + \lVert \omega_s \rVert_{L^{\infty}_x} + \lVert u_n - u_s \rVert^2_{L^{\infty}_x} \right] 
\end{align}

For any time interval $[0,T]$ where $T<\infty$, if we can show that the $L^{\infty}_t H^m_x$ norms of $u_n$ and $u_s$ are bounded, we may substitute these bounds back into \eqref{normal higher-order estimate} and \eqref{super higher-order estimate} and derive upper bounds for $\nu_n \lVert \nabla u_n \rVert^2_{H^m_x}$ and $\nu_s \lVert \nabla u_s \rVert^2_{H^m_x}$, just as was done to arrive at the bounds in \eqref{enstrophy and palinstrophy bounds}. This will show that the $X_n$ and $X_s$ norms of the respective velocity fields are bounded over any time interval $[0,T]$ for every $T<\infty$, allowing us to conclude global solutions of high-regularity. \\

In this regard, we recall a crucial potential theory estimate from \cite{Beale1984RemarksEquations}. \\

\begin{lem} [BKM potential theory estimate] \label{BKM potential theory estimate}
    For $m > \frac{d}{2}+1 = 2$,
    \begin{equation} \label{BKM potential estimate}
        \lVert \nabla u \rVert_{L^{\infty}_x} \lesssim \lVert \omega \rVert_{L^2_x} + \lVert \omega \rVert_{L^{\infty}_x} \left( 1 + \log(1 + \lVert u \rVert_{H^m_x}) \right)
    \end{equation}
\end{lem}

\bigskip

Thus, the quantity in the square brackets in \eqref{higher-order estimate} can be bounded above by:

\begin{align*}
    [\dots] \lesssim (\lVert \omega_n \rVert_{L^2_x} &+ \lVert \omega_s \rVert_{L^2_x}) + \lVert u_n - u_s \rVert^2_{L^{\infty}_x} \nonumber \\
    &+ (\lVert \omega_n \rVert_{L^{\infty}_x} + \lVert \omega_s \rVert_{L^{\infty}_x}) \left( 1+ \log (1 + \lVert u_n \rVert_{H^m_x} + \lVert u_s \rVert_{H^m_x}) \right)
\end{align*}

From \eqref{enstrophy and palinstrophy bounds}, we know that $\lVert \omega \rVert_{L^2_x}$ is bounded above for both fluids. So, we can absorb these upper bounds into the implied multiplicative constant in the upper bound. We can also adjust the argument of the logarithm to make the entire expression in \eqref{higher-order estimate} more amenable to the Gr\"onwall's inequality. For this, we note that for $x\ge 0$,

\begin{equation*}
    \log (1+x) \le \log (1+x^2 + 2x) \le \log \left( 2(1+x^2)\right) \lesssim 1 + \log(1+x^2)
\end{equation*}

From the above arguments,

\begin{align} \label{BKM upper bound simplified}
    [\dots] \lesssim (1 + \lVert u_n - u_s \rVert^2_{L^{\infty}_x} + \lVert \omega_n \rVert_{L^{\infty}_x} + \lVert \omega_s \rVert_{L^{\infty}_x})\left(1 + \log (1 + \lVert u_n \rVert^2_{H^m_x} + \lVert u_s \rVert^2_{H^m_x})\right) 
\end{align}

Denoting $X = 1 + \lVert u_n \rVert^2_{H^m_x} + \lVert u_s \rVert^2_{H^m_x}$, we see that \eqref{higher-order estimate} simplifies to

\begin{equation} \label{higher order skeletal}
    \frac{dX}{dt} \lesssim X (1 + \log X )\left[ 1 + \lVert u_n - u_s \rVert^2_{L^{\infty}_x} + \lVert \omega_n \rVert_{L^{\infty}_x} + \lVert \omega_s \rVert_{L^{\infty}_x} \right]
\end{equation}

At this stage, one can easily apply Gr\"onwall's inequality to draw the conclusion that the $C^0_t H^m_x$ norm is bounded for any finite time $T$, if the quantity in the square brackets in \eqref{higher order skeletal} is integrable in time over $[0,T]$. If that is the case, the required upper bound for $X$ is given by

\begin{equation} \label{final Gronwall's solution}
    \sup_{0\le t \le T} X(t) \le e^{\left(1+\log X(0) \right) e^{\int_0^T[\dots] dt}}
\end{equation}

\subsection{Analogue of the BKM condition} \label{analogue of the BKM condition}

We will now prove the time-integrability condition required to proceed from \eqref{higher order skeletal} to \eqref{final Gronwall's solution}, which will lead us to global solutions from high-regularity data. The claim to be verified is:

\begin{lem} [BKM analogue] \label{BKM analogue}
    For every $0<T<\infty$, 
    \begin{equation}
        \int_0^T \left[ 1 + \lVert u_n - u_s \rVert^2_{L^{\infty}_x} + \lVert \omega_n \rVert_{L^{\infty}_x} + \lVert \omega_s \rVert_{L^{\infty}_x}\right] dt < \infty
    \end{equation}
\end{lem}

    \textit{Proof}:
    \begin{enumerate}
        \item The first term in the integral is obviously finite for every finite $T$. \\
        
        \item For the second term, we use the Gagliardo-Nirenberg interpolation inequality, followed by the Calderon-Zygmund inequality, and the energy boundedness from \eqref{energy and dissipation bound}:
        
        \begin{align*}
            \lVert u_n - u_s \rVert^2_{L^{\infty}_x} &\lesssim \lVert u_n - u_s \rVert_{L^2_x} \lVert D^2 u_n - D^2 u_s \rVert_{L^2_x} \\
            &\lesssim \left( \lVert u_n \rVert_{L^2_x} + \lVert u_s \rVert_{L^2_x} \right) \left (\lVert \nabla \omega_n \rVert_{L^2_x} + \lVert \nabla \omega_s \rVert_{L^2_x} \right) \\
            &\lesssim \left (\lVert \nabla \omega_n \rVert_{L^2_x} + \lVert \nabla \omega_s \rVert_{L^2_x} \right)
        \end{align*}
        
        Thus, integrating over $[0,T]$,
        
        \begin{equation}
            \lVert u_n - u_s \rVert^2_{L^2_t L^{\infty}_x} \lesssim \left( \lVert \nabla \omega_n \rVert_{L^2_t L^2_x} + \lVert \nabla \omega_s \rVert_{L^2_t L^2_x} \right) T^{\frac{1}{2}}
        \end{equation}
        
        This is clearly finite, from the bounds in \eqref{enstrophy and palinstrophy bounds}. \\
        
        \item For the third term, consider the vorticity equation \eqref{normal fluid vorticity}, reproduced here for convenience.
        
        \begin{equation*}
            \partial_t \omega_n + u_n\cdot\nabla\omega_n = \nu_n \Delta \omega_n + \frac{\rho_s}{\rho}\mathcal{T}
        \end{equation*}
        
        Along the characteristics of the flow (denoted by $x^{\alpha}_n (t)$, where $\alpha$ is the initial point), the Duhamel solution to this equation can be written as:
        
        \begin{equation} \label{normal vorticity Duhamel solution}
            \omega_n (t,x^{\alpha}_n(t)) = e^{\nu_n t\Delta} \omega^0_n (\alpha) + \frac{\rho_s}{\rho} \int_0^t e^{\nu_n (t-t')\Delta}  \mathcal{T}(t') \ dt'
        \end{equation}
        
        Now, we make use of an interesting result, the proof of which can be found in Proposition 44 of \cite{Tao254ANew}. \\
        
        \begin{lem} \label{Tao notes lemma}
            In two spatial dimensions, if a field evolves according to the forced heat equation $\partial_t u = \nu\Delta u + F$ (with initial condition $u^0$), then its solution satisfies
            
            \begin{equation}
                \lVert u \rVert_{L^2_t L^{\infty}_x} \lesssim_{\nu} \lVert u^0 \rVert_{L^2_x} + \lVert F \rVert_{L^1_t L^2_x}
            \end{equation}
        \end{lem}
        
        \medskip
        
        Applying this lemma to \eqref{normal vorticity Duhamel solution}, and absorbing the viscosity factor into the implied constant,
        
        \begin{align}
            \lVert \omega_n \rVert_{L^2_t L^{\infty}_x} & \lesssim  \lVert \omega^0_n \rVert_{L^2_x} + \lVert |\omega_s|(\omega_n - \omega_s) \rVert_{L^1_t L^2_x} \nonumber \\
            &\qquad \qquad \qquad \qquad \qquad + \lVert (u_n - u_s)\cdot(B\nabla^{\perp}|\omega_s| - B'\nabla\omega_s) \rVert_{L^1_t L^2_x} \nonumber \\
            &\lesssim  \lVert \omega^0_n \rVert_{L^2_x} + \lVert |\omega_s|(\omega_n - \omega_s) \rVert_{L^1_t L^2_x}
            + \lVert u_n - u_s \rVert_{L^2_t L^{\infty}_x} \lVert \nabla\omega_s \rVert_{L^2_t L^2_x}
        \end{align}
        
        In the above equation, the first term is finite by the high-regularity assumption, and the finiteness of the last term is inferred from the previous analysis combined with the bound from \eqref{enstrophy and palinstrophy bounds}. For the middle term, from H\"older's and Cauchy's inequalities, it is easy to see that
        
        \begin{equation} \label{middle term}
            \lVert |\omega_s|(\omega_n - \omega_s) \rVert_{L^1_t L^2_x} \lesssim \lVert \omega_n \rVert^2_{L^2_t L^4_x} + \lVert \omega_s \rVert^2_{L^2_t L^4_x}
        \end{equation}
        
        Both the terms on the RHS of \eqref{middle term} are handled in the same way, as shown below for the first of the two. From the Ladyzhenskaya inequality in 2D, followed by the Calderon-Zygmund inequality, 
        
        \begin{align*}
            \lVert \omega_n \rVert^2_{L^4_x} &\lesssim \lVert \omega_n \rVert_{L^2_x} \lVert \nabla\omega_n \rVert_{L^2_x} \lesssim \lVert \nabla u_n \rVert_{L^2_x} \lVert \nabla\omega_n \rVert_{L^2_x} \\
            \Rightarrow \lVert \omega_n \rVert^2_{L^2_t L^4_x} &\lesssim \lVert \nabla u_n \rVert_{L^2_t L^2_x} \lVert \nabla\omega_n \rVert_{L^2_t L^2_x} < \infty
        \end{align*}
        
        Once again, from \eqref{energy and dissipation bound} and \eqref{enstrophy and palinstrophy bounds}, we have that the quantities in the RHS are bounded above. Since the third term in Lemma \ref{BKM analogue} can be bounded by $\lVert \omega_n \rVert_{L^2_t L^{\infty}_x} T^{\frac{1}{2}}$, we see that its contribution to the integral is finite. \\
        
        \item That the fourth (and last) term in Lemma \ref{BKM analogue} is finite is proven exactly as the third term above.
    \end{enumerate}
\qed

This completes the proof of the BKM analogue, showing global well-posedness for high-regularity initial data.

\section{Proof of Corollary \ref{GWP smooth data}} \label{GWP for smooth data section}

In the previous sections, we showed that the HVBK equations are globally well-posed for high-regularity data. Now, we seek to show that starting from $C^{\infty}_c$ data, the solutions are globally well-posed and are smooth in space and time. \\

Since a $C^{\infty}_c$ function belongs to $H^m_x$ for every $m$, and data in $H^m_x$ implies global well-posedness, it is easy to see that $C^{\infty}_c$ data means there is a unique global solution in every $H^m_x$, which means it is spatially $H^{\infty}$-smooth, i.e., $u_n,u_s \in C^0([0,\infty[;H^{\infty}(\R^2))$. Time regularity can be shown just as in Section \ref{instantaneous smoothing}. Thus, $u_n,u_s \in C^{\infty}([0,\infty[;H^{\infty}(\R^2))$, implying smoothness of the pressure fields (from elliptic regularity). The uniqueness of the pressure fields can be argued as done in section \ref{regularity of pressure}.

\qed

\appendixtitleon
\begin{appendices}
    \section{Conditions on the vorticity field for finite kinetic energy in 2D incompressible fluids}
    
    It is well-known that an incompressible fluid in $\R^2$ with compactly supported vorticity has finite kinetic energy if and only if the integral of the vorticity vanishes (see Prop. 3.3 in \cite{Majda2002VorticityFlow}). This is due to the slow decay of the Biot-Savart kernel in two dimensions. 
    
    \begin{gather}
        \nabla\times u = \omega \Rightarrow -\Delta u = \nabla\times\omega \nonumber \\
        \therefore u(x) = \frac{1}{2\pi} \int_{\R^2} \frac{(x-y)^{\perp}}{|x-y|^2} \  \omega(y) \ dy \label{velocity from vorticity}
    \end{gather}
    
    However, it turns out that a compactly supported vorticity is not the only way to tackle this problem. In this appendix, we wish to replace this assumption with two weaker conditions: finite enstrophy and finite $L^1$ norm of the first moment of the vorticity.

\begin{lem} \label{finite kinetic energy}
    For an incompressible fluid in 2D, let there be a vorticity field $\omega:\R^2\mapsto\R$ such that $\int_{\R^2} \omega \ dx = 0$, and $\lVert \omega \rVert_{L^2(\R^2)}$ and $\int_{\R^2}\lvert x \rvert \lvert \omega \rvert \ dx$ are both finite. The associated velocity field $u:\R^2\mapsto\R^2$ is defined by \eqref{velocity from vorticity}. Then the kinetic energy of the fluid is finite, i.e., $\lVert u \rVert_{L^2(\R^2)} < \infty$. In particular, $\lVert u \rVert_{L^2(\R^2)} < \lVert \omega \rVert_{L^2(\R^2)}^{\frac{1}{2}} \lVert x \omega \rVert_{L^1(\R^2)}^{\frac{1}{2}}$.
\end{lem}

    \textit{Proof}: Writing \eqref{velocity from vorticity} in Fourier space,
    
    \begin{equation*}
        \hat{u}(k) = \frac{ik^{\perp}}{|k|^2}\hat{\omega}(k)
    \end{equation*}
    
    From Plancherel's theorem, we have
    \begin{equation*}
        \lVert u \rVert_{L^2_x} \lesssim \lVert \hat{u} \rVert_{L^2_k} \lesssim \left( \int \frac{\lvert \hat{\omega}(k) \rvert^2}{|k|^2} dk \right)^{\frac{1}{2}}
    \end{equation*}
    
    We split the integral into one over low frequencies and another over high frequencies. The (to-be-determined) cutoff is denoted $K$.
    
    \begin{equation*}
        \lVert u \rVert_{L^2_x} \lesssim \left( \int_{|k|<K} \frac{\lvert \hat{\omega}(k) \rvert^2}{|k|^2} dk \right)^{\frac{1}{2}} + \left( \int_{|k|\ge K} \frac{\lvert \hat{\omega}(k) \rvert^2}{|k|^2} dk \right)^{\frac{1}{2}}
    \end{equation*}
    
    The high-frequency component is easily seen to be bounded by $\frac{1}{K}\lVert \omega \rVert_{L^2_x}$. The vanishing integral of the vorticity translates to $\hat{\omega}(0) = 0$. Expanding in a Taylor series about $k=0$, for some $0\le \theta \le 1$:
    
    \begin{equation*}
            \hat{\omega}(k) = \hat{\omega}(0) + D\hat{\omega}(\theta k)\cdot k \Rightarrow |\hat{\omega}(k)| \le |D\hat{\omega}(\theta k)| \ |k|
    \end{equation*}

    Also,
    
    \begin{equation*}
        D\hat{\omega}(k_0) = \int i x \omega(x) e^{ik_0 x} dx \Rightarrow \lVert D\hat{\omega} \rVert_{L^{\infty}_k} \le \lVert x \omega \rVert_{L^1_x}
    \end{equation*}
    
    We can thus bound the low-frequency component by $K \lVert x \omega \rVert_{L^1}$. This gives us:
    
    \begin{equation*}
        \lVert u \rVert_{L^2_x} \lesssim K \lVert x \omega \rVert_{L^1_x} + \frac{1}{K} \lVert \omega \rVert_{L^2_x}
    \end{equation*}
    
    Selecting $K =  \lVert x \omega \rVert_{L^1_x}^{-\frac{1}{2}} \lVert \omega \rVert_{L^2_x}^{\frac{1}{2}}$ gives the required result.
        
    \qed
\end{appendices}


\bibliographystyle{spmpsci}
\bibliography{Manuscript}

%
%

%
%

\end{document}